\documentclass[11pt]{amsart}
\usepackage{amssymb}
\usepackage{amsmath}

\newcommand{\bbC}{{\mathbb C}}
\newcommand{\bbZ}{{\mathbb Z}}

\newcommand{\bbP}{{\mathbb P}}

\newcommand{\Oh}{{\mathcal O}}

\DeclareMathOperator{\HH}{H}

\DeclareMathOperator{\II}{\mathcal I_{C/X}}
\DeclareMathOperator{\I^2}{\mathcal I^2_{C/X}}

\DeclareMathOperator{\VV}{V}
\DeclareMathOperator{\UU}{U}
\DeclareMathOperator{\CH}{CH}
\DeclareMathOperator{\Pic}{Pic}

\DeclareMathOperator{\Ext}{Ext}
\DeclareMathOperator{\Hom}{Hom}

\DeclareMathOperator{\im}{Image}

\DeclareMathOperator{\pf}{pf}

\newcommand{\onto}{\twoheadrightarrow}

\newcommand{\into}{\hookrightarrow}
\newcommand{\by}[1]{\xrightarrow{#1}}
\newcommand{\ext}[1]{\operatorname{\stackrel{#1}{\wedge}}}
\newcommand{\tensor}{\otimes}
\newcommand{\isom}{\cong}

\newcommand{\sI}{{\mathcal I_{C/\bbP^4}}}
\newcommand{\sIW}{{\mathcal I_{W/\bbP^4}}}

\newcommand{\sH}{{\mathcal H}}
\newcommand{\sHom}{\mbox{${\sH}{om}$}}

\theoremstyle{plain}
\newtheorem{lemma}{Lemma}
\newtheorem{prop}{Proposition}
\newtheorem{thm}{Theorem}
\newtheorem{conj}{Conjecture}
\newtheorem{cor}{Corollary}

\theoremstyle{remark}
\newtheorem{remark}{Remark}

\newtheorem{claim}{Claim}
\newenvironment{diagram}[1]{\arraycolsep=\doublerulesep\begin{array}{#1}
    }{\end{array}}

\long\def\comment#1{}
\usepackage{MnSymbol}

\begin{document}

\title{Curves on threefolds and a conjecture of
Griffiths-Harris}
\author{G. V. Ravindra}
\address{Department of Mathematics, Indian Institute of
Science, Bangalore -- 560 012, INDIA.}
\email{\tt ravindra@math.iisc.ernet.in}

\begin{abstract}
We prove that any arithmetically Gorenstein curve on a smooth, general
hypersurface $X\subset \bbP^{4}$ of degree at least $6$, is a complete
intersection. This gives a characterisation of complete intersection
curves on general type hypersurfaces in $\bbP^4$. We also verify that
certain $1$-cycles on a general quintic hypersurface are non-trivial
elements of the Griffiths group.
\end{abstract}
\date{\today} 
\maketitle

\section{Introduction}
We work over $\bbC$, the field of complex numbers. By a {\it general}
point of a variety, we shall mean a point in a Zariski open subset and
by a {\it very general} point we mean a point in the complement of a
countable union of proper closed subvarieties.

For a very general hypersurface $X\subset\bbP^3$ of degree at least $4$,
the Noether-Lefschetz theorem (NLT) says that every curve $C\subset X$
is a complete intersection of $X$ with a surface in $\bbP^3$ i.e.,
$C=X\cap S$ where $S\subset \bbP^3$ is a surface. Motivated by this,
Griffiths and Harris \cite{GH1} conjectured that the following
analogue of NLT holds for curves in threefolds.

\begin{conj}[Griffiths-Harris, \cite{GH1}]
Let $X\subset\bbP^4$ be a very general hypersurface of degree $d\geq
6$. Then any curve $C\subset X$ is of the form $C=X\cap S$, where
$S$ is a surface in $\bbP^4$.
\end{conj}

 For the sake of brevity, we shall call curves $C \subset X$ which are
not intersections of $X$ with any surface as {\it special}. Voisin
(see \cite{V}) showed that a general threefold $X\subset\bbP^4$ always
contains special curves $C \subset X$, thus proving that this
conjecture is false. In fact, one can consider the analogous question
for codimension two subvarieties in higher dimensional hypersurfaces;
in \cite{MPR3}, it is shown that there exists a large class of special
codimension two subvarieties in smooth hypersurfaces of dimension at
least three and degree at least two.

The aim of this note is to show that though NLT for curves in surfaces
does not generalise to curves in threefolds, a restricted version of
this theorem related to the non existence of certain special curves on
very general hypersurfaces in $\bbP^3$ also holds for general
hypersurfaces in $\bbP^4$. We shall make this precise now.

We start with a few definitions. A vector bundle $E$ on $X$ is said to
be {\it arithmetically Cohen-Macaulay} (ACM for short) if $\HH^i(X,
E(\nu))=0$, $\forall \,\nu\in\bbZ$ and $ 0<i<\dim{X}$.  Similarly, a
subscheme $Y\subset X$ with ideal sheaf $\mathcal{I}_{Y/X}$ is said to
be ACM if $\HH^i(X,\mathcal{I}_{Y/X}(\nu))=0$, $\forall \,\nu\in\bbZ$
and $1\leq i \leq \dim{Y}$. In addition, if $Y$ has codimension two in
$X$, we shall say $Y$ is {\it arithmetically Gorenstein} if $Y$ is the
zero scheme of a section of a rank two bundle $E$ on $X$. It is not
hard to see in this case that if $X$ is a smooth projective hypersurface
of dimension at least $3$, then $Y$ is a complete intersection if and
only if $E$ is a sum of line bundles and that $Y$ is ACM if and only
if $E$ is ACM.

An equivalent formulation of NLT says that if $X\subset \bbP^3$ is a
very general hypersurface of degree at least $4$, then any line bundle
$L$ on $X$ is $\Oh_X(m)$ for some $m\in \bbZ$; hence $L$ is
ACM. Rephrasing this, we may say that as a consequence of this
theorem, any ACM line bundle on such an $X$ is the restriction of a
line bundle on $\bbP^3$.

One might now wonder that though the analogue of NLT for curves in
threefolds $X\subset\bbP^4$ is false, is it still true that any ACM
rank two vector bundle on $X$ is the restriction of a rank two bundle
on $\bbP^4$. In fact, one might even be tempted to formulate the
``Noether-Lefschetz question'' for higher rank ACM bundles as follows:
Given an ACM rank $r$ bundle $E$ on a very general smooth hypersurface
$X\subset\bbP^n$, is $E$ the restriction of a vector bundle on
$\bbP^n$ ? The case $(r,n)=(1,3)$ is implied by the NLT. It is easy to
see that any such extension, if it exists, is necessarily ACM on
$\bbP^n$.  By a theorem of Horrocks (see \cite{Ho}), any ACM vector
bundle on $\bbP^n$ is a sum of line bundles. Thus the
Noether-Lefschetz question for higher rank ACM vector bundles can also
be viewed as an extension of Horrocks' splitting criterion to bundles
on hypersurfaces $X\subset \bbP^n$.  Notice that the converse, namely
that any sum of line bundles on $X$ extends to $\bbP^n$ for $n\geq 4$
follows by the Grothendieck-Lefschetz theorem.

Buchweitz-Greuel-Schreyer have shown (see \cite{BGS}) that there do
exist non-trivial ACM bundles of sufficiently high rank on any
hypersurface $X$.  Conjecture B of {\it op.~cit.}  tells us precisely
beyond what rank one might expect to get non-trivial ACM bundles.  The
main result of this paper is the following which can be viewed as a
verification of the first non-trivial case of (a strengthening of)
this conjecture.

\begin{thm}\label{mgh}
Let $X$ be a general hypersurface in $\bbP^4$ of degree $d\geq 6$. Any
arithmetically Gorenstein curve $C\subset X$ is a complete
intersection. Equivalently, any ACM bundle $E$ of rank two on $X$ is a
sum of line bundles.
\end{thm}

By remark 3.1 in \cite{CM2}, it then follows that the above result is
true for a general hypersurface of degree $d\geq 6$ in $\bbP^n$ for
$n\geq 4$. However in {\it op. cit.}, it has been shown that the
result is also true for $d=3,4$ and $5$ when $n\geq 5$. Thus we
recover the following theorem proved by us using completely different
methods:

\begin{cor}[Mohan Kumar-Rao-Ravindra, \cite{MPR1}]\label{acminp5}
Any ACM bundle of rank two on a general hypersurface $X\subset\bbP^n$,
$n\geq 5$, of degree at least $3$ is a sum of line bundles.
\end{cor}

Soon after the main steps in this paper were carried out, we were able
to extend the methods of {\it loc.~cit.} to prove theorem \ref{mgh}
(see \cite{MPR2}). Partial results (see \cite{MPR1} for details) in
this direction were also obtained by Chiantini and Madonna
\cite{CM,CM2, CM3}. 

Theorem \ref{mgh} is sharp: A smooth hypersurface in $\mathbb{P}^4$ of
degree $\leq 5$ contains a line. The corresponding rank two bundle,
via Serre's construction, is ACM but not decomposable. Similarly,
there are smooth hypersurfaces in $\mathbb{P}^4$ of degree $\geq 6$
which contain a line. Hence the hypothesis of generality cannot be
dropped.

We briefly outline the proof of Theorem \ref{mgh}. We may assume that
the vector bundle $E$ is normalised (i.e. $h^0(E(-1))= 0$ and
$h^0(E)\neq 0$). By Proposition \ref{reduction}, its first Chern class
$\alpha:=c_1(E)$ satisfies the inequality $3-d\leq \alpha \leq d-2$.
Suppose on the contrary, a general hypersurface supports such a bundle
which is indecomposable. This implies the following: Let $S'$ be an
open set of the parameter space of smooth hypersurfaces of degree $d$
in $\bbP^4$ which support such a bundle and $\mathcal{X}'\to S'$ be
the universal hypersurface. Then there exists a family of rank two
vector bundles $\mathcal{E}\to \mathcal{X}'$ such that $\forall s\in
S'$, $E_s:=\mathcal{E}_{|X_s}$ is a normalised, indecomposable ACM
bundle of rank two on $X_s$ with $c_1(E_s)=\alpha$. Associated to
this, there is a family of null-homologous $1$-cycles $\mathcal{Z}\to
S'$ whose fibre at any point $s\in S'$ is $\mathcal{Z}_s=dC_s-lD_s$
where $C_s\subset X_s$ is the zero locus of a section of $E_s$,
$l=l(s)=\deg{C_s}$ and $D_s\subset X_s$ is a plane section. To such a
family of cycles, one can associate a normal function
$\nu_{\mathcal{Z}}$ and its infinitesimal invariant
$\delta\nu_{\mathcal{Z}}$ (see section \ref{maindefs} for
definitions). By a result of Mark Green \cite{MG} and Voisin
(unpublished), $\delta\nu_{\mathcal{Z}}\equiv 0$ whenever $d\geq
6$. On the other hand, by refining a method of X.~Wu (see \cite{Wu1}),
we show that in the situation described above,
$\delta\nu_{\mathcal{Z}}\not\equiv 0$ when $d\geq 5$. This is a
contradiction when $d\geq 6$.

On hypersurfaces of degree $d\leq 5$, it is easy to construct
(normalised) indecomposable ACM bundles of rank two (see
\cite{BBR}). The (refined) criterion of Wu has an interesting
consequence for such bundles when $d=5$. Recall that the Griffiths
group of codimension $k$ cycles is defined to be the group of
homologically trivial codimension $k$ cycles modulo the subgroup of
cycles algebraically equivalent to zero. As a consequence of the
non-degeneracy of the infinitesimal invariant, we have

\begin{cor}\label{gg}
Let $E$ be a normalised, indecomposable ACM bundle of rank two on
$X_5\subset\bbP^4$, a smooth, general quintic hypersurface. If
$Z=5C-lD$ is as above, then $Z$ defines a non-trivial element of
$\mbox{Griff}^2(X)$, the Griffiths group of codimension $2$ cycles on
$X$.
\end{cor}

The proof of the above corollary is identical to Griffiths' proof (see
\cite{G}) where he shows that the difference of two distinct lines
defines a non-trivial element in the Griffiths group. Hence we shall
only sketch the proof and refer the reader to {\it op.~cit.} for
details. Since $\delta\nu_{\mathcal{Z}}\neq 0$, this implies that
$\nu_{\mathcal{Z}}$ is not locally constant (see \cite{MG}), hence $Z$
has non-trivial Abel-Jacobi image. Now the subgroup of cycles
algebraically equivalent to zero is contained in the kernel of the
Abel-Jacobi map. Hence the corollary.

\subsection*{Comparison with Wu's results:} 
In \cite{Wu}, using his criterion for the non-degeneracy of the
infinitesimal invariant, Wu is able to prove the following:
\begin{thm}
Let $X\subset\bbP^4$ be a general, smooth hypersurface of degree
$d\geq 6$, and let $C\subset X$ be a smooth curve with $\deg{C}\leq
2d-1$. Then $C=X\cap\bbP^2$ is a plane section.
\end{thm}

Thus Theorem \ref{mgh} may be viewed as a generalisation of this
theorem of Wu. Though any characterisation of complete intersection
curves cannot obviously have a constraint on their degrees as in the
above theorem, it is interesting to note that the proof of Theorem
\ref{mgh} also follows by reducing to the case of bounded degree
curves. To see this, let $E$ be a rank two ACM bundle on $X$ and
assume that it has a non-zero section whose zero locus $C$ is a
curve. The Grothendieck-Riemann-Roch formula expresses $\chi(E)$ as a
function of $c_1(E)$ and $c_2(E)$ (see \cite{CM} for a precise
formula). Since $E$ is ACM,
$\chi(E(b))=h^0(E(b))-h^3(E(b))=h^0(E(b))-h^0(E(-c_1+d-5-b))$. Choosing
$b>0$ so that $\chi(E(b))\geq 0$, we see that $c_2(E(b)) $ is bounded by a
function of $c_1(E(b))$. From the outline of the proof given above, we
may assume that $c_1(E)$ is bounded; hence it follows that
$\deg{C}=c_2(E)$ is bounded in terms of the degree of $X$.
\comment{\frac{dc_1^3}{6}+ \frac{(5-d)dc_1^2}{4}
  -\frac{(5-d)c_1}{2}-\frac{c_1c_2}{2} +\frac{dc_1(2d^2-15d+35)}{12}+
  \frac{d(-d^3+10d^2-35d+50)}{12}}

\begin{remark}
Theorem \ref{mgh} has now been generalised to complete intersections
subvarieties of sufficiently high multi-degree in projective space
(see \cite{BR}).
\end{remark}

\section{Preliminaries}
\subsection{Reductions}
Let $X\subset\bbP^4$ be a smooth hypersurface of degree $d$. By the
Grothendieck-Lefschetz theorem, we have $\Pic(X)\isom\bbZ$ with
generator $\Oh_X(1)$. Also by the weak Lefschetz theorem, we have
$\HH^{2i}(X,\bbZ)\isom \HH^{2i}(\bbP^4,\bbZ)\isom\bbZ$ for $i=1,~2.$
With these identifications, we may treat the first and second Chern
classes of any vector bundle $E$ on $X$ as integers.

In this section, we shall show that it is enough to consider rank two
ACM bundles whose first Chern class $\alpha$ satisfies the inequality
$3-d \leq \alpha\leq d-2$. A useful result that we shall use is the
following remark which can be found in \cite{Ha2}.

\begin{lemma}\label{minimal}
Let $E$ be a normalised, indecomposable ACM bundle of rank two on a
smooth projective variety $X$ with $\Pic(X)\isom\bbZ$. Then the zero
scheme of any non-zero section of $E$ has codimension $2$ in $X$. In
particular, if $\dim{X}\geq 2$, the zero scheme is non-empty.
\end{lemma}

Let $X\subset \bbP^4$ be a smooth hypersurface of degree $d$, $E$ an
ACM bundle of rank two on $X$. If $C$ is the zero scheme of a section
of $E$ (hence a curve by Lemma \ref{minimal}), one has a short exact
sequence 
\begin{equation}\label{serre}
0 \to \Oh_X \to E \to \II(\alpha) \to 0,
\end{equation}
where $\mathcal{I}_{C/X}$ denotes the ideal sheaf of $C$ in $X$. It
follows that such a $C$ is sub-canonical i.e., $\omega_C\isom\Oh_C(m)$
for some $m\in\bbZ$. To determine $m$, note that from the above exact
sequence, $E(-\alpha)\tensor\Oh_C\isom\II/\I^2$. Taking determinants
on both sides and using adjunction, we have
$K_X\tensor\omega_C^{-1}=\det{E}\tensor\Oh_C(-2\alpha)$. Rewriting, we
get $m=\alpha +d -5$.

The inclusions $C\subset X \subset \bbP^4$ yield the following
short exact sequences:
\begin{equation}\label{curveinp4}
 0 \to \sI \to \Oh_{\bbP^4} \to \Oh_C \to 0,
\end{equation}
\begin{equation}\label{curveinx}
 0 \to \II \to \Oh_{X} \to \Oh_C \to 0,
\end{equation}
\begin{equation}\label{inclofideals}
0 \to \Oh_{\bbP^4}(-d) \to \sI \to \II \to 0.
\end{equation}
Finally, $E$ has a length one resolution by a sum of line bundles on
$\bbP^4$:
\begin{equation}\label{acmres}
0 \to F_1 \by{\Phi} F_0 \to E \to 0
\end{equation}
where $F_0=\bigoplus_{i=1}^r\Oh_{\bbP^4}(-a_i)$, $a_i\geq 0$ for
all $i$, $F_1=F_0^{\vee}(\alpha-d)$, and $\Phi$ is a
skew-symmetric matrix. Details may be found in \cite{Beau, MPR1}.

Recall (see \cite{Mumford}) that a coherent sheaf $\mathcal{F}$ on $X$
is said to be $m$-regular in the sense of Castelnuovo-Mumford if
$\HH^i(X, \mathcal{F}(m-i))=0$ for $i>0$. When $m=0$, we say that
$\mathcal{F}$ is regular.

\begin{lemma}\label{reg}
With notation as above, $E(d-\alpha-1)$ is regular in the sense of
Castelnuovo-Mumford.
\end{lemma}

\begin{proof}
We need to check that $\HH^i(X,E(d-\alpha-1-i))=0$ for $i=1,2,3$.  The
vanishings for $i=1, 2$ follow from the fact that $E$ is ACM.  For
$i=3$, note that $\HH^3(E(d-\alpha-1-3)) \isom \HH^0(E^{\vee}(\alpha
-d +4 +d -5)) \isom \HH^0(E(-1))=0$ where the first isomorphism is by
Serre duality and the second follows from the fact $E^{\vee} \isom
E(-\alpha)$.
\end{proof}

\begin{prop}\label{reduction}
Let $E$ be a normalised, indecomposable ACM bundle of rank two on
$X\subset\bbP^4$, a general hypersurface of degree $d$ at least
$6$. Then its first Chern class satisfies the inequality, $3-d \leq
\alpha \leq d-2$.
\end{prop}

\begin{proof}
$E(d-\alpha -1)$ is regular implies that it is globally generated (see
  \cite{Mumford}, page 99). Since $h^0(E(-1))=0$, we must have
  $d-\alpha -1 > -1$ and so $\alpha < d$. When $\alpha=d-1$, $E$ is
  regular and so all its (minimal) generators are in degree $0$ and we
  have a resolution
$$0 \to F_1=\Oh_{\bbP^4}(-1)^{2d} \by{\Phi} F_0=\Oh_{\bbP^4}^{2d} \to
E \to 0.$$
This implies in particular that $X$ is defined by $\pf(\Phi)=0$ where
for any skew-symmetric matrix $M$, $\pf(M):=\sqrt{\det M}$. By an easy
dimension count (or see Corollary 2.4 in \cite{Beau}), a general
hypersurface of degree at least $6$ is not a linear Pfaffian and hence
$X$ cannot support such an $E$.

To see the lower bound, we reproduce the argument from
\cite{MPR1}. Let $C$ be a the zero-scheme of a section of $E$ and let
$\pi:C \to \bbP^1$ be a general projection so that $\pi$ is
finite. Then
$$\pi_{\ast}\omega_C \isom \sHom(\pi_{\ast}\Oh_C, \Oh_{\bbP^1}(-2)).$$
Since $\Oh_{\bbP^1} \subset \pi_{\ast}\Oh_C$ is a direct summand, this
implies that $\HH^0(\pi_{\ast}\omega_C(2))\isom\HH^0(\omega_C(2))\neq
0$.  On the other hand, since $C$ is ACM, it is clear from the
cohomology sequence associated to sequence (\ref{curveinx}) that
$\HH^0(\Oh_C(l))=0$ if $l <0$. Putting these together, we get 
$\alpha+ d -5 +2 \geq 0$ or equivalently $\alpha \geq 3-d$.
\end{proof}

\begin{cor}\label{correg}
With notation as above,
\begin{itemize}
\item[a)] $F_0=\bigoplus_{i=1}^r\Oh_{\bbP^4}(-a_i)$, $0 \leq a_i < d
-\alpha$.
\item[b)] $F_1\isom\bigoplus_i\Oh_{\bbP^4}(-b_i), \hspace{3mm}
  b_i=d-a_i-\alpha > 0$.
\item[c)] $\HH^0(F_1)=0$ and $\HH^0(F_0) \isom \HH^0(E)$.
\end{itemize}
\end{cor}

\begin{proof}
These facts follow immediately from the regularity of $E(d-\alpha -1)$
and the fact that $F_1\isom F_0^{\vee}(\alpha -d)$.
\end{proof}

\subsection{Griffiths' infinitesimal invariant and a 
result of Green and Voisin.}\label{maindefs}

The main object of study is an invariant defined by Griffiths
\cite{G1, G2} which we briefly discuss now. For the details, we refer
the reader to {\it op.~cit.} or Chapter 7 of Voisin's book
\cite{Vbook}.

Let $\mathcal{X} \to S$ be the universal family of smooth, degree $d$
hypersurfaces in $\bbP^{2m}$. Let $\mathcal{C}\subset\mathcal{X}$ be a
family of codimension $m$ subvarieties over $S$. For a point $s\in S$,
let $X=X_{s}$ and $C:=\mathcal{C}_{s} \subset X$. If $l$ is the degree
of $C$ and $D_s$ is a codimension $m$ linear section, then the family
of cycles $\mathcal{Z}$ with fibre
$\mathcal{Z}_s:=d~\mathcal{C}_s-lD_s$ for $s\in S$ defines a
(fibre-wise null-homologous) cycle in
$\CH^m(\mathcal{X}/S)_{hom}$. Let $\mathcal{J}:=\{J(X_s)\}_{s\in S}$
be the family of intermediate Jacobians. In such a situation,
Griffiths defines a holomorphic function $\nu_{\mathcal{Z}}:S \to
\mathcal{J}$, called the {\it normal function}, by
$\nu_{\mathcal{Z}}(s)=\mu_s(Z_s)$ where
$\mu_s:\CH^m(X_s)_{\mbox{\small{hom}}} \to J(X_s)$ is the {\it
  Abel-Jacobi} map from the group of null-homologous cycles to the
intermediate Jacobian. This normal function satisfies a
``quasi-horizontal'' condition (see \cite{Vbook}, Definition 7.4).
Associated to the normal function $\nu_{\mathcal{Z}}$ above, Griffiths
(see \cite{G1} or \cite{Vbook} Definition 7.8) has defined the
infinitesimal invariant $\delta\nu_{\mathcal{Z}}$. Later Green
\cite{MG} generalised this definition and showed that Griffiths'
original infinitesimal invariant is just one of the many infinitesimal
invariants that one can associate to a normal function. He also showed
that in particular $\delta\nu_{\mathcal{Z}}(s)$ is an element of the
dual of the middle cohomology of the following complex
$$ \ext{2}\HH^1(X,T_X)\tensor\HH^{m+1,m-2}(X)\to
  \HH^1(X,T_X)\tensor\HH^{m,m-1}(X)\to
  \HH^{m-1,m}(X).
$$

We now specialise to the case $m=2$ where $X\subset\bbP^4$ is a smooth
hypersurface and $C\subset X$ is a curve of degree $l$. Then
$Z:=dC-lD$ is a nullhomologous $1$-cycle with support $W:=C\bigcup D$.
At a point $s\in S$, this infinitesimal invariant is a functional
$$\delta\nu_{\mathcal{Z}}(s):
\ker\left(\HH^1(X,T_X)\tensor\HH^1(X,\Omega^2_X) \to
\HH^2(X,\Omega^1_X)\right)\to\bbC.$$ The following result of Griffiths gives
an explicit formula for computing the infinitesimal invariant
associated to the family $\mathcal{Z}$ at a point $s\in S$ when
restricted to
$$\ker\left(\HH^1(X,T_X)\tensor\HH^1(X,\mathcal{I}_{W/X}\tensor\Omega^2_X) \to
\HH^2(X,\Omega^1_X)\right).$$

\begin{thm}[Griffiths \cite{G1,G2}]
Let $\nu_{\mathcal{Z}}$ be the normal function as described
above. Consider the following diagram:
\begin{equation}\label{formula}
\begin{diagram}{ccccccccc}
& & & & \HH^1(X, T_X)\tensor\HH^1(X,\mathcal{I}_{W/X}\tensor\Omega_X^2) & &
& &
\\
& & & & \downarrow{\beta} & \searrow{\gamma}& & & \\
0 & \to & \HH^1(W,\Omega^1_X\tensor\Oh_W)/\HH^1(X,\Omega^1_X)
&
\by{\lambda} &\HH^2(X,\mathcal{I}_{W/X}\tensor\Omega^1_X) & \to &
\HH^2(X,\Omega^1_X) &
\to & 0 \\
& & \downarrow{\chi} & & & &  & & \\
& & \bbC & & & & & & \\
\end{diagram}
\end{equation}
where $\chi$ is given by integration over the cycle $Z$. Then
$\delta\nu_{\mathcal{Z}}(s)$, the infinitesimal invariant evaluated at
a point $s\in S$, is the composite
$$ \ker{\gamma}\to \frac{\HH^1(W,\Omega^1_X\tensor\Oh_W)}
{\HH^1(X, \Omega^1_X)}\by{\chi} \bbC,$$
where the map $$ \ker{\gamma}\to \frac{\HH^1(W,
  \Omega^1_X\tensor\Oh_W)} {\HH^1(X, \Omega^1_X)}$$ is induced by the
map $\beta$ and the above short exact sequence.
\end{thm}
The map $\chi$ can be understood as follows. Since $D$ is a general
plane section of $X$, by Bertini $C\cap{D}=\emptyset$.  Thus
$\Oh_W\isom\Oh_C\oplus\Oh_D$ and so
$$\HH^1(W, \Omega^1_X\tensor\Oh_W)\isom \HH^1(C,
\Omega^1_X\tensor\Oh_C) \oplus \HH^1(D, \Omega^1_X\tensor\Oh_D). $$
For any irreducible curve $T\subset X$, let
$$r_T:\HH^1(T, \Omega^1_X\tensor\Oh_T) \to
\HH^1(T, \Omega^1_T)\isom \bbC$$ denote the natural restriction map.
 For any element $(a,b)\in
\HH^1(W, \Omega^1_X\tensor\Oh_W)$, we define
\begin{equation}\label{chi}
\chi(a,b):=dr_C(a)-lr_D(b) \in \bbC.
\end{equation} It is clear that this map
factors via the quotient
$\HH^1(W, \Omega^1_X\tensor\Oh_W)/\HH^1(X, \Omega^1_X).$

The following result is due to Green \cite{MG} and Voisin
(unpublished).

\begin{thm}
Let $X\subset\bbP^4$ be a general hypersurface of degree at least
$6$. Then the infinitesimal invariant $\delta\nu$ of any
quasi-horizontal normal function $\nu$, is zero.
\end{thm}

\subsection{Wu's criterion}\label{setup}
Now we are ready to prove the final step i.e. that there are no
non-trivial normalised ACM bundles $E$ of rank two on a general
hypersurface $X\subset\bbP^4$ of degree $d\geq 6$ such that $3-d\leq
\alpha \leq d-2$. We shall suppose the contrary: that such an $E$
exists on a general hypersurface $X$ as above. In such a situation,
(see section 3 of \cite{MPR1} for details) there exists a rank two
bundle $\mathcal{E}$ on the universal hypersurface $\mathcal{X}
\subset \bbP^4 \times S'$ where $S'$ is a Zariski open subset of $S$,
the moduli space of smooth, degree $d$ hypersurfaces of $\bbP^4$, such
that for a general point $s\in S'$, $\mathcal{E}_{|X_s}$ is
normalised, indecomposable, ACM with first Chern class
$\alpha$. Furthermore, from the construction of this family, one sees
that there exists a family of curves $\mathcal{C} \to S'$ such that
$\mathcal{C}_s$ is the zero locus of a section of
$\mathcal{E}_{|X_s}$. Let $\mathcal{Z}$ be a family of $1$-cycles with
fibre $\mathcal{Z}_s:=d\mathcal{C}_s-lD_s$ where, as before, $D_s$ is
plane section of $X_s$ and $l=l(s)$ is the degree of $\mathcal{C}_s$.

We shall show that under the hypotheses above,
$\delta\nu_{\mathcal{Z}} \not\equiv 0$. The non-degeneracy of the
infinitesimal invariant is shown by refining Xian Wu's proof in
\cite{Wu1}. The proof has three main steps, which we describe now.

Let $\partial_f:\Omega^3_{\bbP^4}(2d)\to K_{\bbP^4}(3d)$ be the
exterior differential between sheaves of meromorphic differential
forms, where $\Omega_{\mathbb{P}^4}^3(2d)$ is identified with the
sheaf of meromorphic $3$-forms with poles of order at most $2$ along
$X$ and $K_{\bbP^4}(3d)$ is identified with the sheaf of meromorphic
$4$-forms with poles of order at most $3$ along $X$. Composing with
the natural map $ K_{\bbP^4}(3d)\onto K_{\bbP^4}(3d)/ K_{\bbP^4}(2d)$,
we get a map $\bar\partial_f:\Omega^3_{\bbP^4}(2d)\to
K_{\bbP^4}(3d)/K_{\bbP^4}(2d)$. Using the identification
$\Omega^3_{\bbP^4}\isom T_{\bbP^4}\tensor K_{\bbP^4}$, and taking
cohomology, we get
$$\HH^0(\bbP^4, T_{\bbP^4}\tensor K_{\bbP^4}(2d))\by{\partial_f}
\HH^0(\bbP^4, K_{\bbP^4}(3d)) \to \frac{\HH^0(\bbP^4,
  K_{\bbP^4}(3d))}{\HH^0(\bbP^4, K_{\bbP^4}(2d))}.$$ The cokernel of
the composite map above can be identified with $\HH^2(X, \Omega^1_X)$
(see \cite{CGGH}, Page 174 or \cite{JL}, Chapter 9 for details). Let
$\bar{U}\subset\HH^0(\bbP^4, K_{\bbP^4}(3d))$ be the subspace defined
as follows:
\begin{equation}\label{ubar}
\bar{U}:=\partial_f\HH^0(\bbP^4, T_{\bbP^4}\tensor K_{\bbP^4}(2d))\cap
\HH^0(\bbP^4, \mathcal{I}_{W/\bbP^4}\tensor K_{\bbP^4}(3d)).
\end{equation}

The key ingredient in the proof is the following commutative
diagram (see {\it op.~ cit.}):
\begin{equation}
\begin{array}{ccc}
\HH^0(\bbP^4, \Oh_{\bbP^4}(d))\tensor\HH^0(\bbP^4,\, \mathcal
I_{W/\bbP^4}\tensor K_{\bbP^4}(2d)) & \stackrel{\gamma'}\to &
\frac{\HH^0(\bbP^4, K_{\bbP^4}(3d))}{\partial_f\HH^0(\bbP^4,
  T_{\bbP^4}\tensor K_{\bbP^4}(2d))}\\ \downarrow & & \downarrow
\\ \HH^1(X, T_{X})\tensor \HH^1(X,\,
\mathcal{I}_{W/X}\tensor\Omega^2_{X}) & \by{\gamma} & \HH^2(X,
\Omega^1_{X}).\\
\end{array}
\end{equation}
Here the right vertical map is as explained above. The horizontal maps
$\gamma$ and $\gamma'$ are (essentially) cup product maps. The
vertical map on the left is a tensor product of two maps.  The first
factor is the composite $\HH^0(\bbP^4, \Oh_{\bbP^4}(d))\to \HH^0(X,
\Oh_X(d))\to\HH^1(X, T_X)$. The normal bundle of $X\subset \bbP^4$ is
$\Oh_X(d)$ and $\HH^0(X, \Oh_X(d))\to \HH^1(X, T_X)$ is the natural
coboundary map in the cohomology sequence of the tangent bundle
sequence for this inclusion. The second factor is the composite
$$\HH^0(\bbP^4,\, \mathcal I_{W/\bbP^4}\tensor K_{\bbP^4}(2d))\to
\HH^0(X,\, \mathcal I_{W/\bbP^4}\tensor K_{\bbP^4}(d)\tensor\Oh_X)\to
\HH^1(X,\, \mathcal{I}_{W/\bbP^4}\tensor K_{\bbP^4}(d)\tensor{T}_X).$$
Here the first map is the natural restriction map and the second is
obtained as above by first tensoring the tangent bundle sequence with
$K_{\bbP^4}(d)\tensor\mathcal{I}_{W/\bbP^4}$ and observing that (i)
$\mathcal{I}_{W/\bbP^4}\tensor\Oh_X\isom\mathcal{I}_{W/X}$ and,
(ii) $T_X\tensor K_{\bbP^4}(d)\isom \Omega^2_X$.  \comment{ Tensor
  with, the tangent bundle sequence for the inclusion
  $X\subset\bbP^4$. By, we get a coboundary map
$$\HH^0(X, \mathcal I_{W/\bbP^4}\tensor K_{\bbP^4}(d)\tensor\Oh_X) \to \HH^1(X,
\mathcal{I}_{W/\bbP^4}\tensor\Omega^2_{X}) = \HH^1(X,
\mathcal{I}_{W/X}\tensor\Omega^2_{X}).$$ }

This diagram yields a map $\ker\gamma' \to \ker\gamma$. To show that
the infinitesimal invariant $\delta\nu_{\mathcal{Z}}(s): \ker\gamma
\to \bbC$ is non-zero, we shall show that the composite map
\begin{equation}\label{keysurj}
\ker\gamma' \to \ker\gamma \to \bbC 
\end{equation}
is non-zero (= surjective).

This is done as follows: consider the exact sequence
$$0 \to \Oh_X(-d) \to \Omega^1_{\bbP^4{|X}} \to \Omega^1_X \to 0.$$
Taking second exterior powers and tensoring the resulting sequence by
$\Oh_X(d)$, we get a short exact sequence
\begin{equation}\label{2ndext}
0 \to \Omega^1_X \to \Omega^2_{\bbP^4{|X}}(d)\to\Omega^2_X(d)\to 0.
\end{equation}
The inclusion $\Omega^1_{X|W} \into \Omega^2_{\bbP^4}(d)_{|W}$ induces
a map of cohomologies and we let
$$\VV:=\ker[\HH^1(W,\Omega^1_{X{|W}}) \to
\HH^1(W,\Omega^2_{\bbP^4}(d)_{|W})].$$

The surjectivity of the composite map in equation (\ref{keysurj}) in
turn is accomplished by constructing a surjection from $\ker\gamma'$
to the vector space $\bar U$ (defined in equation \eqref{ubar}) such
that this fits into a commutative diagram
\begin{equation}
\begin{array}{ccccc}
\ker\gamma' & \to & \ker\gamma & & \\
\twoheaddownarrow & & &  \searrow & \\
\bar{U} & \onto & V & \onto & \bbC \\
\end{array}
\end{equation}
where the map $\ker\gamma \to \bbC$ is the infinitesimal invariant
evaluated at the point $s\in S$.  In the next section, we shall carry
out the three steps viz,
\begin{enumerate}
\item[Step 1.] There exists a surjection $\chi:V\onto\bbC$.
\item[Step 2.] There exists a surjection $\bar{U}\onto \VV$.
\item[Step 3.] There exists a surjection $\ker\gamma'\onto\bar{U}$.
\end{enumerate}

\section{Proof of Theorem \ref{mgh}}
\subsection{Step 1: The surjection $\chi:V\onto\bbC$}

The main result of this section is the following

\comment{
To show that
$$\chi:\VV \to \bbC$$ is surjective, we shall need to state and prove
a few lemmas which we do now.
We now prove the following key lemma.}

\begin{prop}\label{keymaplemma}
The restriction map
$$\HH^1(X, \Omega^2_{\bbP^4}(d){|_{X}}) \to \HH^1(C,
\Omega^2_{\bbP^4}(d){|_{C}})$$ is zero.
\end{prop}

Let us see how this Proposition implies Step 1. The natural map
$\Omega^1_{X{|C}} \to \Omega^1_C$ yields a push-out diagram for
sequence (\ref{2ndext}) (see \cite{Rotman} pages 41--43 for a
definition):
\begin{equation}\label{pushout}
\begin{diagram}{ccccccccc}
0 & \to & \Omega^1_{X{|C}} & \to & \Omega^2_{\bbP^4}(d)_{|C} &
\to &
\Omega^2_X(d)_{|C} & \to & 0 \\
& & \downarrow & & \downarrow & & || & & \\
0 & \to & \Omega^1_C & \to &
\mathcal{F} & \to & \Omega^2_X(d)_{|C} & \to & 0. \\
\end{diagram}
\end{equation}

\begin{lemma}
The map $\HH^1(C,\Omega^1_C) \to \HH^1(C,\mathcal{F})$ in the
associated cohomology sequence of the bottom row in diagram
(\ref{pushout}) is zero. Thus we have a surjection
$$\VV_C:=\ker[\HH^1(C,\Omega^1_{X_{|C}}) \to
\HH^1(C,\Omega^2_{\bbP^4}(d)_{|C})] \onto \HH^1(C,\Omega^1_C).$$
\end{lemma}

\begin{proof}
We have a commutative diagram
\[\begin{array}{ccc}
\HH^1(X,\Omega^1_X) & \to & \HH^1(X,\Omega^2_{\bbP^4}(d)_{|_X}) \\
\downarrow &  & \downarrow \\
\HH^1(C,\Omega^1_{X_{|C}}) & \to &
\HH^1(C,\Omega^2_{\bbP^4}(d)_{|_C}) \\
\downarrow &  & \downarrow \\
\HH^1(C,\Omega^1_C) & \to & \HH^1(C,\mathcal{F}). \\
\end{array}\]
The composite of the vertical maps on the left is the natural
restriction map $\HH^1(X,\Omega^1_X) \to \HH^1(C,\Omega^1_C)$ which
maps $h_X\mapsto h_C$ where $h_A$ is the hyperplane class for any
scheme $A$. Since both these cohomologies are one-dimensional with
$h_X$ and $h_C$ as the respective generators, this map is an
isomorphism. Now $\HH^1(X,\Omega^2_{\bbP^4}(d)_{|_X}) \to
\HH^1(C,\Omega^2_{\bbP^4}(d)_{|_C})$ is the zero map by Proposition
\ref{keymaplemma}, and so this implies that the map
$\HH^1(C,\Omega^1_C) \to \HH^1(C,\mathcal{F})$ is zero. Thus we have a
surjection
$
\VV_C \onto \HH^1(C,\Omega^1_C)$.
\end{proof}

\begin{cor}[Step 1]
The composite map 
\[
\begin{array}{ccc}
\VV_C \into \VV=\ker[\HH^1(W,\Omega^1_{X_{|W}}) \to
\HH^1(W,\Omega^2_{\bbP^4}(d)_{|W})] & \stackrel{\chi}\to & \bbC \\
\end{array}
\]
is a surjection. Hence $\chi$ is a surjection.
\end{cor}

\begin{proof}
This first inclusion follows from the fact that $\Oh_W \isom \Oh_C
\oplus \Oh_D$. The surjectivity of the composite follows from the
definition of $\chi$ (see equation (\ref{chi})) and the above lemma.
\end{proof}

To prove Proposition \ref{keymaplemma}, we shall need a few more
results which we prove now.

Applying the functor $\sHom_{\Oh_{\bbP^4}}(-,\Oh_{\bbP^4})$ to
sequence (\ref{acmres}), we get (see \cite{MPR1})
\begin{eqnarray}\label{dualres}
0 \to F_0^{\vee} \by{\Psi} F_1^{\vee} \to E^{\vee}(d) \to 0 .
\end{eqnarray}

Let $\phi: F_0^{\vee} \to \Oh_{\bbP^4}$ be any morphism (equivalently
a section $\phi \in \HH^0(F_0)$). Associated to any such morphism, we
have a push-out diagram:

\[
\begin{array}{ccccccccc}
0 & \to & F_0^{\vee} & \by{\Psi} &  F_1^{\vee} & \to & E^{\vee}(d)
&
\to & 0 \\
& & \downarrow{\phi} & & \downarrow & & || & & \\
0 & \to & \Oh_{\bbP^4} & \to &  G & \to &  E^{\vee}(d) &
\to & 0.\\
\end{array}
\]

Conversely, we have

\begin{lemma}\label{comm}
Any exact sequence
$$0  \to  \Oh_{\bbP^4} \to   G  \to   E^{\vee}(d) \to 0,$$
arises as a push-out diagram above.
\end{lemma}

\begin{proof}
Any exact sequence as above corresponds to an element of
$\Ext^1_{\Oh_{\bbP^4}}(E^{\vee}(d), \Oh_{\bbP^4})$. Applying the
functor $\Hom_{\Oh_{\bbP^4}}(-,\Oh_{\bbP^4})$ to sequence
(\ref{dualres}) we get,
$$0 \to \Hom_{\Oh_{\bbP^4}}(F_1^{\vee},\Oh_{\bbP^4}) \to
\Hom_{\Oh_{\bbP^4}}(F_0^{\vee},\Oh_{\bbP^4}) \to
\Ext^1_{\Oh_{\bbP^4}}(E^{\vee}(d), \Oh_{\bbP^4}) \to 0.$$
The first term is $\HH^0(F_1)$ which is zero by Corollary
\ref{correg} and thus we have an isomorphism $\HH^0(F_0)
\isom\Ext^1_{\Oh_{\bbP^4}}(E^{\vee}(d),\Oh_{\bbP^4})$.
\end{proof}

Putting these together, we get the following
\begin{cor}
Let $E^{\vee} \by{s^{\vee}} \Oh_X$ be the map induced by a section
$s\in \HH^0(E)$. Consider the pull-back diagram (see \cite{Rotman}
pages 51-53 for definition)
\[
\begin{array}{ccccccccc}
0 & \to & \Oh_{\bbP^4} & \to &  G & \to &  E^{\vee}(d) & \to & 0 \\
& & ||  & & \downarrow & & \downarrow{s^{\vee}} & & \\
0 & \to & \Oh_{\bbP^4} & \to & \Oh_{\bbP^4}(d) & \to & \Oh_X(d) &
\to & 0. \\
\end{array}
\]

By Lemma \ref{comm}, there is a section $\phi\in \HH^0(F_0)$ such
that the following diagram commutes:
\begin{equation}
\begin{array}{ccccccccc}\label{keycomm}
0 & \to & F_0^{\vee} & \by{\Psi} &  F_1^{\vee} & \to & E^{\vee}(d)
&
\to & 0 \\
& & \downarrow{\phi}  & & \downarrow & & \downarrow{s^{\vee}} & & \\
0 & \to & \Oh_{\bbP^4} & \to & \Oh_{\bbP^4}(d) & \to & \Oh_X(d) &
\to & 0. \\
\end{array}
\end{equation}
In fact, under the isomorphism $\HH^0(\bbP^4, F_0) \isom \HH^0(X,E)$
(Corollary \ref{correg} (c)), $\phi$ maps to $s$.
\end{cor}

\begin{remark}\label{splitsurj}
Since $F_0^\vee=\bigoplus\Oh_{\bbP^4}(a_i)$ where $a_i\geq 0$, the map
$\phi$ restricted to a summand $\Oh_{\bbP^4}(a_i)$ with $a_i>0$ is
zero. Hence $\phi$ is a split surjection.
\end{remark}

\begin{proof}[Proof of Proposition \ref{keymaplemma}]
We remark that $\HH^1(C, \Omega^2_{\bbP^4}(d){|_{C}})=0$ when $\alpha
< 2$, and so the lemma is obvious in these cases. The proof for all
values $\alpha < d-1$ is as follows. From Corollary \ref{correg},
$F_1^\vee=\bigoplus_i\Oh_{\bbP^4}(b_i)$, $b_i>0$ and so by Bott's
formula (see for instance, \cite{OSS} Page 8) $\HH^i(\bbP^4,
\Omega^2_{\bbP^4}\tensor F_1^{\vee})=0$ for $i=1,~2$. This implies
that the boundary map $\HH^1(X, \Omega^2_{\bbP^4}\tensor E^{\vee}(d))
\to \HH^2(\bbP^4, \Omega^2_{\bbP^4}\tensor F_0^{\vee})$ in the
cohomology sequence associated to sequence (\ref{dualres})
$\tensor\Omega^2_{\bbP^4}$ is an isomorphism.

Next, we tensor diagram (\ref{keycomm}) by $\Omega^2_{\bbP^4}$ and
take cohomology to get a commutative diagram
\[
\begin{array}{ccc}
\HH^1(X, \Omega^2_{\bbP^4}\tensor E^{\vee}(d)) & \isom &
\HH^2(\bbP^4, \Omega^2_{\bbP^4}\tensor F_0^{\vee})\\
\downarrow & & \downarrow \\ 
\HH^1(X, \Omega^2_{\bbP^4}\tensor\Oh_{X}(d)) & \isom &
\HH^2(\bbP^4, \Omega^2_{\bbP^4}),\\
\end{array}
\]
where the isomorphism in the bottom row follows again from Bott's
formula ({\it op.~cit.}).  By Remark \ref{splitsurj}, the right
vertical map above is onto, and this implies that the map
\begin{equation}\label{etoox}
\HH^1(X, \Omega^2_{\bbP^4}\tensor E^{\vee}(d)) \to
\HH^1(X, \Omega^2_{\bbP^4}(d){|_{X}})
\end{equation}
is onto. The map $E^\vee(d)\to \Oh_X(d)$ in diagram (\ref{keycomm}) is
induced by a section $s\in\HH^0(X,E)$ and hence has image
$\mathcal{I}_{C/X}(d)$, where $C=Z(s)$. Thus the map in equation
(\ref{etoox}) factors via $\HH^1(X,
\II\tensor\Omega^2_{\bbP^4}(d)_{|X})$ and so the map
$$\HH^1(X, \II\tensor\Omega^2_{\bbP^4}(d)_{|X}) \to
\HH^1(X, \Omega^2_{\bbP^4}(d){|_{X}})$$ is surjective. Thus we are done.
\end{proof}

\begin{remark}
Proposition \ref{keymaplemma} is a crucial refinement of Wu's
criterion. Wu actually requires that $\HH^1(C,
\Omega^2_{\bbP^4}(d)_{|C})=0$. Since $\HH^1(X,
\Omega^2_{\bbP^4}(d)_{|X})$ is $1$-dimensional, we were hopeful that
the weaker statement that the map is zero, which is really what we
need, might hold with our hypotheses.
\end{remark}

\subsection{ Step 2: The surjection $\bar{U}\to \VV$}
We first describe the map $\bar{U} \to \VV$.

Tensoring the short exact sequence
$$ 0 \to  T_X \to  T_{{\bbP^4|X}} \to \Oh_X(d) \to 0$$ with
$K_{\bbP^4}(2d)_{|W}$ and taking cohomology, we get 
$$0 \to \HH^0(W, T_{X}\tensor K_{\bbP^4}(2d)_{|W}) \to
\HH^0(W, T_{\bbP^4}\tensor K_{\bbP^4}(2d)_{|W}) \to
\HH^0(W, K_{\bbP^4}(3d)_{|W}).$$
Since $ T_X\tensor K_{\bbP^4}(d) \isom \Omega^2_X$, we have
the following commutative diagram:
\begin{equation}\label{unexplained}
\begin{array}{ccccccc}
0 & \to & \UU  & \to & \HH^0(\bbP^4, T_{\bbP^4}\tensor
K_{\bbP^4}(2d))
& \to & \HH^0(W, K_{\bbP^4}(3d)_{|W}) \\
& & \downarrow & & \downarrow & & || \\
0 & \to & \HH^0(W, \Omega^2_{X}(d)_{|W}) & \to &
 \HH^0(W, T_{\bbP^4}\tensor K_{\bbP^4}(2d)_{|W}) & \to &
\HH^0(W, K_{\bbP^4}(3d)_{|W}).\\
\end{array}
\end{equation}
Here $U$ is defined so that the top row is left exact. From the
exactness of the cohomology sequence associated to sequence
(\ref{2ndext}), we get
$$\im[\HH^0(W, \Omega^2_{X}(d)_{|W}) \to \HH^1(W, \Omega^1_{X|W})] =
\ker[\HH^1(W, \Omega^1_{X|W}) \to \HH^1(W,
  \Omega^2_{\bbP^4}(d)_{|W})]=\VV,$$ and hence a surjective map
$\HH^0(W, \Omega^2_{X}(d)_{|W}) \onto \VV.$ Consider the composite
$$\UU \to \HH^0(W, \Omega^2_{X}(d)_{|W}) \onto \VV.$$

\begin{cor}
The map $U\to\VV$ factors as $U \to \bar{U} \to \VV$.
\end{cor}

\begin{proof}
Let $\widetilde{U}$ be the kernel of the map $\HH^0(\bbP^4,
T_{\bbP^4}\tensor K_{\bbP^4}(2d)) \to \HH^0(X,
K_{\bbP^4}(3d)_{|X})$. Looking at the diagram analogous to
(\ref{unexplained}) obtained by replacing $W$ by $X$, we see that
there is a map $\widetilde{U}\to \HH^0(X, \Omega^2_X(d))$. The
boundary map $\HH^0(X, \Omega^2_X(d))\to \HH^1(X, \Omega^1_X)$ in the
cohomology sequence associated to diagram (\ref{2ndext}) is the zero
map (this is because the composite map $\HH^1(X, \Omega^1_X) \to
\HH^1(X, \Omega^2_{\bbP^4|X}(d))\isom \HH^2(\bbP^4,
\Omega^2_{\bbP^4})$ is the Gysin isomorphism). This implies that the
map $U \to V$ above factors as $U \onto U/\widetilde{U}\to V$.  Next
we claim that the map $U \onto U/\widetilde{U}$ factors as $U \to
\bar{U}\to U/\widetilde{U}$. For this we define
$U_P:=\ker[\HH^0(\bbP^4, T_{\bbP^4}\tensor K_{\bbP^4}(2d))\by{\partial_f}
  \HH^0(\bbP^4, K_{\bbP^4}(3d))]$. It is enough to check the following:
\begin{enumerate}
\item $U_P\subset \widetilde{U} \subset U$.
\item $\partial_f$ restricts to a surjective map $U \to \bar{U}$ which
  induces an isomorphism $U/U_P\isom \bar{U}$.
\end{enumerate}
These follow easily from the definitions of $U_P$, $\bar{U}$ and
$\tilde{U}$.
\end{proof}

To show that the map $\bar{U}\to\VV$ defined above is a surjection, it
is enough to prove that the map $U \to \VV$ is surjective. For this,
we shall need some vanishings which we prove now.

\begin{lemma}\label{c4}
For $d\geq 3$, $\HH^1(\bbP^4,\mathcal
I_{W/\bbP^4}(2d-4))=0=\HH^2(\bbP^4,\mathcal I_{W/\bbP^4}(2d-5))$.
\end{lemma}

\begin{proof}
Tensoring the exact sequence
\begin{equation}\label{resforD1}
0 \to \Oh_X(-2) \to \Oh_X(-1)^{\oplus 2} \to \mathcal I_{D/X} \to 0 \, ,
\end{equation}
by $\II$, we get the exact sequence
$$0\to \II(-2) \to \II(-1)^{\oplus 2} \to \mathcal I_{W/X} \to 0 .$$
Left exactness here can be checked at the level of stalks using the
fact that $C\cap D=\emptyset$. 

For the first vanishing, since $C$ is ACM, taking cohomology of the
above sequence we see that the boundary map
$\HH^1(X,\mathcal{I}_{W/X}(2d-4))\to \HH^2(X,\mathcal I_{C/X}(2d-6))$
in the cohomology sequence associated to the above exact sequence is
an injection. Using the exact sequence $0 \to \mathcal{I}_{C/X} \to
\Oh_X \to \Oh_C \to 0$, we see that $\HH^2(X,\mathcal
I_{C/X}(2d-6))\isom \HH^1(C,\Oh_C(2d-6))$ which in turn is Serre dual
to $\HH^0(C,\Oh_C(\alpha-d+1))$.  \comment{We claim that this group is
  zero. The claim will follow from the following fact: {\it for any
    ACM curve $Y\subset X$, $\HH^0(Y,\Oh_Y(k))=0$ for $k<0$.}  To see
  this consider the sequence
$$0 \to \mathcal I_{Y/X}(k) \to \Oh_X(k) \to \Oh_Y(k) \to 0.$$
Taking cohomology, we get a sequence
$$\cdots \to \HH^0(X,\Oh_X(k)) \to \HH^0(Y,\Oh_Y(k)) \to
\HH^1(X,\mathcal I_{Y/X}(k)) \to \cdots.$$ The first term is clearly
zero and the second is zero since $Y$ is ACM. Thus the middle term
vanishes.}  Since $C$ is ACM and $\alpha < d-1$, we have
$\HH^0(C,\Oh_C(\alpha-d+1))=0$, and so
$0=\HH^1(X,\mathcal{I}_{W/X}(2d-4))=\HH^1(\bbP^4,\mathcal
I_{W/\bbP^4}(2d-4))$ (for the last equality, use sequence
(\ref{inclofideals}) with $C$ replaced by $W$).

For the second vanishing, consider the short exact sequence
$$ 0 \to \mathcal{I}_{W/\bbP^4} \to \Oh_{\bbP^4} \to \Oh_W \to 0.$$
Taking cohomology, it is easy to see that there are isomorphisms
$$\HH^2(\bbP^4,\mathcal I_{W/\bbP^4}(2d-5))\isom
\HH^1(W,\Oh_W(2d-5))\isom\HH^1(C,\Oh_C(2d-5)) \oplus
\HH^1(D,\Oh_{D}(2d-5)).$$ The first term is Serre dual to
$\HH^0(C,\Oh_C(\alpha-d))$ and the second term to
$\HH^0(D,\Oh_D(2-d))$. Since $C$, $D$ are ACM, $\alpha < d-1$ and
$d\geq 3$, it follows that $\HH^0(C,\Oh_C(\alpha-d))$ and
$\HH^0(D,\Oh_D(2-d))$ are both zero. This finishes the proof.
\end{proof}

\begin{lemma}\label{c2}
For $d\geq 5$, $\HH^1(\bbP^4, T_{\bbP^4}\tensor\mathcal
I_{W/\bbP^4}(2d-5))=0.$
\end{lemma}

\begin{proof}
Tensoring the Euler sequence by $\sIW(2d-5)$, we get a short exact
sequence (see \cite{Wu} for left exactness)
$$0  \to  \sIW(2d-5)  \to  \sIW(2d-4)^{\oplus{5}}  \to
\sIW(2d-5)\tensor T_{\bbP^4}  \to  0. $$
This gives rise to a part of a long exact sequence of cohomology
$$ \to \HH^1(\bbP^4, \sIW(2d- 4))^{\oplus{5}} \to \HH^1(\bbP^4,
T_{\bbP^4}\tensor\sIW(2d-5)) \to \HH^2(\bbP^4, \sIW(2d-5)) \to $$ By
Lemma \ref{c4}, the extreme terms vanish and so we are done.
\end{proof}

\begin{prop}\label{utovsurj}
For $d\geq 5$, the natural map $U \to V$ is a surjection.
\end{prop}

\begin{proof}
The middle vertical arrow in diagram (\ref{unexplained}) can be seen
to be a surjection by using the fact that the cokernel of this map
injects into $\HH^1(\bbP^4, T_{\bbP^4}\tensor
K_{\bbP^4}\tensor\mathcal I_{W/\bbP^4}(2d))$ which vanishes by Lemma
\ref{c2}. By snake lemma, the first map is also a surjection. Thus the
map $U \to \HH^0(W, \Omega^2_{X}(d)_{|W})$ is a surjection. This
finishes the proof.
\end{proof}

Thus we have the required
\begin{cor}[Step 2]
For $d\geq 5$, the map $\bar{U}\to\VV$ is also a surjection.
\end{cor}

\subsection{Step 3: The surjection $\ker\gamma' \to \bar{U}$}
We first describe the map $\ker\gamma'\to \bar{U}$.

Recall from section \ref{setup} that $\gamma'$ is the natural map
$$\HH^0(\bbP^4, \Oh_{\bbP^4}(d))\tensor\HH^0(\bbP^4, \mathcal I_{W/\bbP^4}\tensor
K_{\bbP^4}(2d))  \to 
\frac{\HH^0(\bbP^4, K_{\bbP^4}(3d))}{\partial_f\HH^0(\bbP^4, T_{\bbP^4}\tensor
K_{\bbP^4}(2d))}\,\, \cdot$$
Consider the multiplication map
\begin{equation}\label{*}
\HH^0(\bbP^4, \Oh_{\bbP^4}(d))\tensor\HH^0(\bbP^4, \mathcal I_{W/\bbP^4}\tensor
K_{\bbP^4}(2d))\to
\HH^0(\bbP^4, \mathcal{I}_{W/\bbP^4}\tensor{K_{\bbP^4}}(3d)).
\end{equation}
 Restricting
this map to $\ker{\gamma'}$, we get a map
$$\ker{\gamma'}\to \bar{U}=\partial_f\HH^0(\bbP^4, T_{\bbP^4}\tensor
K_{\bbP^4}(2d))\cap \HH^0(\bbP^4, \mathcal{I}_{W/\bbP^4}\tensor
K_{\bbP^4}(3d)).$$ 

\begin{prop}[Step 3]
For $d\geq 5$, the map $\ker\gamma' \to \bar{U}$ is surjective.
\end{prop}

\begin{proof}
To prove surjectivity of the above map, it is enough to prove that for
$d\geq 5$, the multiplication map in equation (\ref{*}) i.e., the map
$$\HH^0(\bbP^4,\Oh_{\bbP^4}(d))\tensor\HH^0(\bbP^4,\mathcal{I}_{W/\bbP^4}(2d-5))
\to \HH^0(\bbP^4,\mathcal{I}_{W/\bbP^4}(3d-5)) $$ is surjective. 

To see this, we first tensor the exact sequence 
\begin{equation}\label{resforD}
0 \to \Oh_X(-2) \to \Oh_X(-1)^{\oplus 2} \to \mathcal{I}_{D/X} \to 0 \, ,
\end{equation}
by $E$ to get
\begin{equation}\label{aa}
0 \to E(-2) \to E(-1)^{\oplus 2} \to \mathcal{I}_{D/X}E \to 0 \, .
\end{equation}

Let $T_m:=\HH^0(X, \Oh_X(m))$. The exact sequence above gives rise to
a diagram with exact rows:
\[
\begin{array}{ccccccccc}
0 & \to & \HH^0(X, E(n-2))\tensor T_m & \to & 
\HH^0(X, E(n-1))^{\oplus{2}}\tensor T_m & \to  &
\HH^0(X, \mathcal I_{D/X}E(n))\tensor T_m & \to & 0 \\
& &   \downarrow  &  & \downarrow &  &  \downarrow  & & \\
0 & \to & \HH^0(X, E(m+n-2)) & \to & 
\HH^0(X, E(m+n-1))^{\oplus{2}} & \to &
\HH^0(X, \mathcal I_{D/X}E(m+n)) & \to & 0. \\
\end{array}
\]
Here the vertical arrows are all multiplication maps and the exactness
on the right is because $E$ is ACM.

Since $E$ is $(d-\alpha -1)$--regular, the middle vertical arrow is a
surjection for $n\geq d-\alpha$ and $m\geq 0$. It follows that the
multiplication map
$$\HH^0(X,\mathcal{I}_{D/X}E(n))\tensor\HH^0(X,\Oh_X(m))\to
\HH^0(X,\mathcal{I}_{D/X}E(m+n))$$ is surjective for $n\geq
(d-\alpha)$ and $m\geq 0$. Next consider the exact sequence
$$0 \to \mathcal{I}_{D/X} \to \mathcal{I}_{D/X}E \to
\mathcal{I}_{W/X}(\alpha) \to 0$$ obtained by tensoring sequence
(\ref{serre}) by $\mathcal{I}_{D/X}$. As before, left exactness here
can be checked at the level of stalks using the fact that
$C\cap{D}=\emptyset$. Repeating the previous argument, it is easy to
check that the multiplication map
$$\HH^0(X,\mathcal{I}_{W/X}(n))\tensor\HH^0(X,\Oh_X(m))\to
\HH^0(X,\mathcal{I}_{W/X}(m+n))$$ is surjective for $n\geq d$ and
$m\geq 0$. In particular, if $d\geq 5$, the map is surjective for
$n=2d-5$. Also the multiplication map
$$\HH^0(\bbP^4, \Oh_{\bbP^4}(m))\tensor\HH^0(\bbP^4, \Oh_{\bbP^4}(n))
\to \HH^0(\bbP^4, \Oh_{\bbP^4}(m+n))$$ is surjective for $m,n\geq 0$.
Now using the exact sequence
$$ 0 \to \Oh_{\bbP^4}(-d) \to \mathcal{I}_{W/\bbP^4} \to
\mathcal{I}_{W/X} \to 0,$$ and repeating the argument above, we can
conclude using snake lemma, that the multiplication map
$$\HH^0(\bbP^4,\mathcal{I}_{W/\bbP^4}(2d-5))\tensor\HH^0(\bbP^4,\Oh_{\bbP^4}(d))
\to \HH^0(\bbP^4, \mathcal{I}_{W/\bbP^4}(3d-5))$$ is surjective (again
$d\geq 5$ is needed here).
\end{proof}

\subsection{Non-degeneracy of the infinitesimal invariant}
\begin{prop}\label{finale}
In the situation above, $\delta\nu_{\mathcal{Z}} \not\equiv 0$.
\end{prop}

\begin{proof}
We shall show that $\delta\nu_{\mathcal{Z}}(s)\neq 0$ at any point
$s\in S$ parametrising a smooth hypersurface $X\subset \bbP^4$.  From
steps 1--3, we have surjections $\ker\gamma'\onto \bar{U} \onto V
\stackrel{\chi}\onto \bbC$. By the compatibility of these maps (see
\cite{Wu1}) with the map $\ker\gamma'\to\ker\gamma$ and those in
diagram (\ref{formula}), we conclude (using Griffiths' formula) that
$\delta\nu_{\mathcal{Z}}(s)\neq 0$.
\end{proof}

\begin{proof}[Proof of Theorem \ref{mgh}]
Assume that a general hypersurface $X$ supports an indecomposable, ACM
vector bundle $E$ of rank two. As seen earlier, we may assume that $E$
is normalised, with $3-d\leq\alpha\leq d-2$. Let $\mathcal{Z}$ be the
family of degree zero $1$-cycles defined earlier. By the refined Wu's
criterion $\delta\nu_{\mathcal{Z}}\not\equiv 0$: this contradicts
Green's theorem. Thus we are done.
\end{proof}
\comment{
\section{Cycles on the quintic threefold}
As mentioned in the introduction, it is easy to construct
indecomposable, normalised  ACM bundles of rank two on a general
hypersurface of degree $d\leq 5$. Let $E$ be such a bundle on
$X\subset\bbP^4$, a general quintic hypersurface and so $c_1(E)\leq
d-2=3$. Let $C$ be the zero locus of a section and $Z:=dC-lD$ be a
homologically trivial cycle. As before, let $\mathcal{Z}$ be the
family of homologically trivial $1$-cycles for which $Z$ is the fibre
at a general point.

\begin{proof}[Proof of Corollary \ref{gg}]
The proof consists of two parts. The first part is to show that the
(rational equivalence class of the) cycle $Z$ is non-zero as an
element of $\CH_1(X)$. The second part is to show that $Z$ is not
algebraically equivalent to zero. For the first part: by the results
of the earlier section, the infinitesimal invariant
$\delta\nu_{\mathcal{Z}}\neq 0$. By \cite{MG}, this means that the
normal function $\nu_{\mathcal{Z}}$ is not locally constant. In
particular, on a general hypersurface $X$, the image of $Z$ under the
Abel-Jacobi map is non zero and so $Z$ determines a non-trivial cycle
in the degree zero part of $\CH_1(X)$. The second part uses a theorem
of Griffiths \cite{G} which states that the kernel of the Abel Jacobi
map contains the subgroup of cycles algebraically equivalent to
zero. Thus $Z$ is not algebraically equivalent to $0$. Since $Z$ is
homologically trivial by definition, it is a non-trivial element of
the Griffiths group.
\end{proof}
}

\section*{Acknowledgements} 
The present work has greatly benefited from our collaboration with
N.~Mohan Kumar and A.~P.~Rao. We are grateful to them for sharing
their ideas, and for the education we received from them during that
time. We thank L.~Chiantini and C.~Madonna from whose papers we learnt
some very useful results. Jishnu Biswas initially, and two unknown
referees later, invested considerable time and effort in a careful
reading of the manuscript. We are grateful to them for their detailed
comments which led to a better understanding of various issues and an
improvement in the exposition.


\begin{thebibliography}{}

\bibitem{BBR} Biswas, I.; Biswas, J.; Ravindra, G.~V., {\it On some
moduli spaces of stable vector bundles on cubic and quartic
threefolds}, Journal of Pure and Applied Algebra 212 (2008), No. 10,
2298--2306.

\bibitem{BR} Biswas, J., Ravindra, G.~V., {\it Arithmetically
Cohen-Macaulay bundles on complete intersection varieties of
sufficiently high multi-degree}, Preprint.

\bibitem{CGGH} Carlson, James; Green, Mark; Griffiths, Phillip;
Harris, Joe, {\it Infinitesimal variations of Hodge structure. I},
Compositio Math.  50 (1983), no. 2-3, 109--205.

\bibitem{Beau} Beauville, Arnaud., {\it Determinantal hypersurfaces},
Dedicated to William Fulton on the occasion of his 60th birthday.
Michigan Math. J.  48 (2000), 39--64.

\bibitem{BGS} Buchweitz, R.-O.; Greuel, G.-M.; Schreyer, F.-O., {\it
Cohen-Macaulay modules on hypersurface singularities. II},
Invent. Math.  88 (1987), no. 1, 165--182.

\bibitem{CM} Chiantini, Luca; Madonna, Carlo., {\it A splitting
criterion for rank 2 bundles on a general sextic threefold},
Internat. J. Math.  15 (2004), no. 4, 341--359.

\bibitem{CM2} Chiantini, L.; Madonna, C. K., {\it ACM bundles on
general hypersurfaces in $\bbP^5$ of low degree}, Collect. Math.
56 (2005), no. 1, 85--96.

\bibitem{CM3} Chiantini, Luca ; Madonna, Carlo, {\it ACM bundles on
a general quintic threefold}, Dedicated to Silvio Greco on the occasion
of his 60th birthday (Catania, 2001).  Matematiche (Catania) 55
(2000), no. 2, 239--258 (2002).

\bibitem{MG} Green, Mark L., {\it Griffiths' infinitesimal invariant
and the Abel-Jacobi map}, J. Differential Geom. 29 (1989), no. 3,
545--555.

\comment{ 
\bibitem{G-MS} Green, Mark; M\"uller-Stach, Stefan, {\it Algebraic
cycles on a general complete intersection of high multi-degree of a
smooth projective variety}, Compositio Math.  100 (1996), no. 3,
305--309.  }

\bibitem{G} Griffiths, Phillip A., {\it On the periods of certain
rational integrals I, II},  Ann. of Math. (2) 90 (1969), 460-495;
ibid. (2) 90 1969 496--541.

\bibitem{G1}Griffiths, Phillip A., {\it Infinitesimal variations of
Hodge structure. III. Determinantal varieties and the
infinitesimal invariant of normal functions}, Compositio Math.  50
(1983), no. 2-3, 267--324.

\bibitem{G2} Griffiths, Phillip A., {\it Infinitesimal invariant of
normal functions}, Topics in transcendental algebraic geometry,
305--316, Ann. of Math. Stud., 106, Princeton Univ. Press,
Princeton, NJ, 1984.

\bibitem{GH1} Griffiths, Phillip; Harris, Joe, {\it On the
Noether-Lefschetz theorem and some remarks on codimension-two
cycles}, Math. Ann. 271 (1985), no. 1, 31--51.

\bibitem{Ho} Horrocks, G.,  {\it Vector bundles on the punctured
spectrum of a local ring}, Proc. London Math. Soc. (3) 14 (1964)
689--713.

\bibitem{Ha2} Hartshorne, Robin., {\it Stable vector bundles of rank
$2$ on $\bbP^{3}$},  Math. Ann.  238 (1978), no. 3, 229--280.

\bibitem{JL} Lewis, James D., {\it A survey of the Hodge conjecture},
Second edition. Appendix B by B. Brent Gordon. CRM Monograph Series,
10. American Mathematical Society, Providence, RI, 1999. xvi+368.

\bibitem{M} Madonna, C., {\it A splitting criterion for rank 2 vector
bundles on hypersurfaces in $\bbP^4$}, Rend. Sem. Mat. Univ.
Politec. Torino 56 (1998), no. 2, 43--54 (2000).

\bibitem{MPR1} Mohan Kumar,~N., Rao, A.~P., and Ravindra, G.~V., {\it
  Arithmetically Cohen-Macaulay bundles on hypersurfaces}, Comment.
  Math Helvetici 82 (2007), no. 4, 829--843.

\bibitem{MPR2} Mohan Kumar,~N., Rao, A.~P., and Ravindra, G.~V., {\it
  Arithmetically Cohen-Macaulay bundles on three dimensional
  hypersurfaces}, International Math Research Notices 2007 no. 8,
  Art. ID rnm025, 11 pp.

\bibitem{MPR3} Mohan Kumar,~N., Rao, A.~P., and Ravindra, G.~V., {\it
  On codimension two subvarieties of hypersurfaces}, Proceedings of
  the Bloch Conference, Fields Institute Communications, eds. Rob de
  Jeu and James Lewis.

\bibitem{Mumford} Mumford, David, {\it Lectures on curves on an
  algebraic surface, With a section by G. M. Bergman,} Annals of
  Mathematics Studies, No. 59 Princeton University Press, Princeton,
  N.J.

\bibitem{OSS} Okonek, Christian; Schneider, Michael; Spindler, Heinz,
  {\it Vector bundles on complex projective spaces}, Progress in
  Mathematics 3, Birkh\"auser.

\bibitem{Rotman} Rotman, Joseph J., {\it An introduction to
homological algebra,} Pure and Applied Mathematics, 85. Academic
Press, New York-London, 1979.


\bibitem{V} Voisin, Claire, {\it Sur une conjecture de Griffiths et
  Harris}, Algebraic curves and projective geometry (Trento, 1988),
  270--275, Lecture Notes in Math., 1389, Springer, Berlin, 1989.

\bibitem{Vbook} Voisin, Claire, {\it Hodge theory and complex
algebraic geometry. II,} Translated from the French by Leila
Schneps. Cambridge Studies in Advanced Mathematics, 77. Cambridge
University Press, Cambridge, 2003.


\bibitem{Wu1} Wu, Xian, {\it On an infinitesimal invariant of normal
functions}, Math. Ann.  288 (1990), no. 1, 121--132.

\bibitem{Wu} Wu, Xian, {\it On a conjecture of Griffiths-Harris
generalizing the Noether-Lefschetz theorem}, Duke Math. J.  60
(1990), no. 2, 465--472.

\end{thebibliography}
\end{document}